\documentclass[12pt]{amsart}

\usepackage[usenames,dvipsnames,svgnames,table]{xcolor}

\usepackage{xargs}
\usepackage[colorlinks=true, pdfstartview=FitV,,, linkcolor=blue, citecolor=blue, urlcolor=blue]{hyperref}
\usepackage{fancyref}
\usepackage[colorinlistoftodos,prependcaption,textsize=tiny,linecolor=red,backgroundcolor=red!25,bordercolor=red]{todonotes}
\usepackage[doublespacing]{setspace}

\newtheorem{theorem}{Theorem}

\theoremstyle{definition}

\newtheorem{Definition}[theorem]{Definition}

\usepackage{graphicx} 
\usepackage{subfigure}
\usepackage{amssymb}
\usepackage{amsmath}
\usepackage{amsthm}
\usepackage{amsfonts}
\usepackage{amscd}
\usepackage{caption}
\usepackage{color}
\usepackage{mathtools}
\usepackage{mathrsfs}
\usepackage[mathscr]{eucal}
\usepackage{tabularx}

\definecolor{mylinkcolor}{RGB}{0,0,140}
\hypersetup{colorlinks,allcolors=mylinkcolor,citecolor=mylinkcolor}
\usepackage[capitalize]{cleveref}
\crefname{section}{Section}{Sections}

\usepackage{vmargin}
\setpapersize{USletter}
\setmargrb{2.5cm}{2cm}{2.5cm}{2.cm}
\usepackage{version}
\excludeversion{sagecode}

\usepackage[inline]{enumitem}

\usepackage{tikz}
\usepackage{tkz-graph}
\usepackage{tkz-berge}
\usetikzlibrary{arrows,shapes}

\usepackage{epsf}

\newcommand{\img}{\mathrm{img}}

\definecolor{cof}{RGB}{219,144,71}
\definecolor{pur}{RGB}{186,146,162}
\definecolor{greeo}{RGB}{91,173,69}
\definecolor{greet}{RGB}{52,111,72}

\definecolor{ly}{rgb}{1.00,0.90,0.71}
\definecolor{dr}{rgb}{0.54,0.25,0.27}
\definecolor{db}{rgb}{0,0.2,0.6}
\definecolor{fg}{rgb}{0,0.5,0}
\definecolor{dg}{rgb}{0,0.8,0.6}
\definecolor{g}{rgb}{0.5,1,0}
\definecolor{cy}{rgb}{0,1,1}
\definecolor{ye}{rgb}{1,1,0}
\definecolor{sb}{rgb}{0,0,0.5}
\definecolor{pk}{rgb}{1,0.75,0.8}
\definecolor{r}{rgb}{1,0,0}
\definecolor{pp}{rgb}{0.5,0,0.5}

\begin{document}
	
	\textbf{title: } Persistent hypergraph homology and its applications 
	
	\textbf{keywords: }Topological data analysis, persistent homology, hypergraphs, face-to-face contacts, population behavior.
	
	\textbf{authors: }
	\begin{itemize}
		\item Yaru Gao, School of Mathematical Sciences, Dalian University of Technology, Dalian, Liaoning, P.R.China, 116024
		
		gaoyaru@dlut.edu.cn
		
		\item Yan Xu, School of Mathematical Sciences, Dalian University of Technology, Dalian, Liaoning, P.R.China, 116024
		
		xuyannn@yeah.net
		
		\item Fengchun Lei, School of Mathematical Sciences, Dalian University of Technology, Dalian, Liaoning, P.R.China, 116024
		
		fclei@dlut.edu.cn
	\end{itemize}
	
	\textbf{corresponding author: }Yaru Gao
	
	\textbf{Abstract: }
	Persistent homology theory is a relatively new but powerful method in data analysis. Using simplicial complexes, classical persistent homology is able to reveal high dimensional geometric structures of datasets, and represent them as persistent barcodes. However, many datasets contain complex systems of multi-way interactions, making these datasets more naturally and faithfully modeled by hypergraphs. In this article, we investigate the persistent hypergraph model, an important generalization of the classical persistent homology on simplicial complexes. We introduce a new homology, $\hat{H}$, on hypergraphs and an efficient algorithm to compute both persistent barcodes and $\hat{H}$ barcodes. As example, our theory is demonstrated by analyzing face-to-face interactions of different populations. The datasets that we select consist of baboons in primate center, people from rural Malawi, scientific conference, workplace and high school. 
	
	\section{Introduction}
	
	Topological data analysis (TDA) is a rapidly developing component in modern data science. Its strength has been demonstrated in a variety of area such as image recognition$^{1}$, feature identification$^{2}$, protein folding prediction$^{3}$, medicine design$^{4-7}$ and much more. Persistent homology, a new branch of algebraic topology, has been proposed to bridge traditional topology and geometry. The essential idea is to introduce a filtration process and measure homology groups by their “lifespans” during the process. Different from traditional topological models, the “lifespan” measurement provides a family of geometric characterizations of the topological invariants. Recently, there has been a burgeoning movement in extending persistent homology theory in various directions to suit broader applications. Notably, Carlsson develops multi-dimensional persistent homology theory$^{8}$, Meehan connects quiver representation with persistent homology$^{9}$, the theory also gets combined with sheaf theory$^{10}$, etc$^{11-15}$.

	This article investigates an important generalization called the hypergraph model, and apply it in analysis of contact patterns. The persistent hypergraph theory is originally introduced by Wu and others, together with theoretical results including its stability$^{16-20}$. In this paper, we introduce a new homology, $\hat{H}$, as a complement to the embedded homology. The $\hat{H}$ homology characterizes those higher dimensional simplex whose boundary does not exist. We also provide an efficient algorithm that computes persistent embedded homology and $\hat{H}$ homology simultaneously for hypergraph model.
	
	Measuring and analyzing interactions between individuals provide key information on social contact behaviors in communities. Recently developed technology of wearable sensors makes the collection of reliable contact data a relatively easy task$^{21-26}$. Previously, contact data is mostly analyzed using graphic models where each vertex represents an individual and each edge represents an interaction. The graphic model efficiently demonstrates 0 and 1 dimensional features, in the topological sense, such as connectivity. However, it does not contain any information from higher dimensions. As a generalization, the hypergraph model consists of connectivity information of all dimensions, that is all group meetings of more than 2 individuals are recorded. Moreover, our algorithm efficiently computes the topological features, namely persistent $H^{\inf}$ and $\hat{H}$ homology, from the persistent hypergraph model. Consequently, we obtain new features of contact patterns that are inaccessible and largely neglected in the past.
	
	\section{Background: Classical persistent homology}
	
	\subsection{From point clouds to complexes.}
	The most obvious way to convert a point cloud $X$ in a metric space into a topological object is to use the point cloud as the vertices of a combinatorial graph whose edges are determined by proximity (vertices within some specified distance $r$). Such a graph, while capturing connectivity data, ignores a wealth of higher-order features beyond clustering. These features can be accurately discerned by thinking of the graph as a scaffold for a higher-dimensional object. Specifically, one completes the graph to a \textbf{simplicial complex} — a space built from simple pieces (simplices) identified combinatorially along faces. The choice of how to fill in the higher-dimensional simplices of the proximity graph allows for different global representations. The most natural method is to generate Rips complexes, defined as follows.
	\begin{Definition}
		A point configuration is a finite subset of $\mathbb{R}^n$. Let $X$ be a point configuration, and let $r>0$, the \textbf{Rips complex}, $R(X,r)$ is the simplicial complex whose $k$-simplices are $\{\{x_0,x_1,\dots,x_k\}:x_i\in X, d(x_i,x_j)<r\}$, where $d$ is the Euclidean metric on $\mathbb{R}^n$. 
	\end{Definition}
	To analyze simplicial complexes, algebraic topology offers a mature set of tools for counting and collating holes and other topological features in spaces and maps between them. In the context of high-dimensional data, algebraic topology works like a telescope, revealing objects and features not visible to the naked eye. In what follows, we concentrate on homology for its balance between ease of computation and topological resolution$^{27}$.
	
	Let $C_*(X,r)$ be the chain complex associated to $R(X,r)$ i.e. we have
		$$C_0\xleftarrow{\partial_1}C_1\xleftarrow{\partial_2}C_2\xleftarrow{\partial_3}\cdots$$ where $C_k$ is the free abelian group generated by $k$-simplices of $R(X,r)$ and $\partial_k$ is the $k$-th boundary map. We denote the cycles by $Z_k(X,r)=\ker(\partial_k)$, the boundaries by $B_k(X,r)=\text{img}(\partial_{k+1})$, the homology by $H_k(X,r)=Z_k(X,r)/B_k(X,r)$ and the number of generators of homology group $H_k$ is the Betti number $\beta_k$.
		
	\subsection{Persistence}
	Despite being both computable and insightful, the homology of a complex associated to a point cloud at a particular $r$ is insufficient: it is a mistake to ask which value of $r$ is optimal. One requires a means of declaring which holes are essential and which can be safely ignored. 
	
	Persistence, as introduced by Edelsbrunner, Letscher, and Zomorodian$^{29}$ and refined by Carlsson and Zomorodian$^{30}$, is a rigorous response to this problem.
	
	\begin{Definition}
		Let $0<r_1<r_2<\cdots$, the \textbf{persistent complex} $\mathscr{C}(X)$ is the set of chain complexes $C_*(X,r_1),C_*(X,r_2),\dots$ together with chain maps $f^i:C_*(X,r_i)\to C_*(X,r_{i+1})$ where $f^i$ consists of the natural embeddings $f^i_k:C_k(X,r_i)\xhookrightarrow[]{} C_k(X,r_{i+1})$.
	\end{Definition}

	To reveal which features persist, we examine the induced inclusion maps $f_*^i:H_*(X,r_i)\to H_*(X,r_{i+1})$.
	\begin{Definition}
		Let $0<r_1<r_2<\cdots$ and $i<j$, the $(i,j)$ \textbf{persistent homology} of $\mathscr{C}(X)$, denoted $H^{i, j}_* (X)$ is defined  to be the image of the induced homomorphism $f_*^{i,j}:H_*(X,r_i)\to H_*(X,r_{j})$. In the other words, $H_k^{i,j}(X)=Z_k(X,r_i)/(Z_k(X,r_i)\cap B_k(X,r_j))$. The number of generators of homology group $H_k^{i,j}$ is called the persistent Betti number $\beta_k^{i,j}$.
	\end{Definition}
	
	We can visualize the persistent homology by barcodes. For example, we have 21 points in $\mathbb{R}^2$. Figure 1 represents persistent homology of the Rips complexes based on these points. The first row represents the 21 points with different radii. In the persistent barcodes representation, the horizontal axis corresponds to radius and the bars represent unordered generators of homology.
	
     \begin{center}
     	\begin{figure}[htbp!]\label{fig:vr}
     		\begin{tikzpicture}
     			\node at (0,0) {\includegraphics[width=80pt]{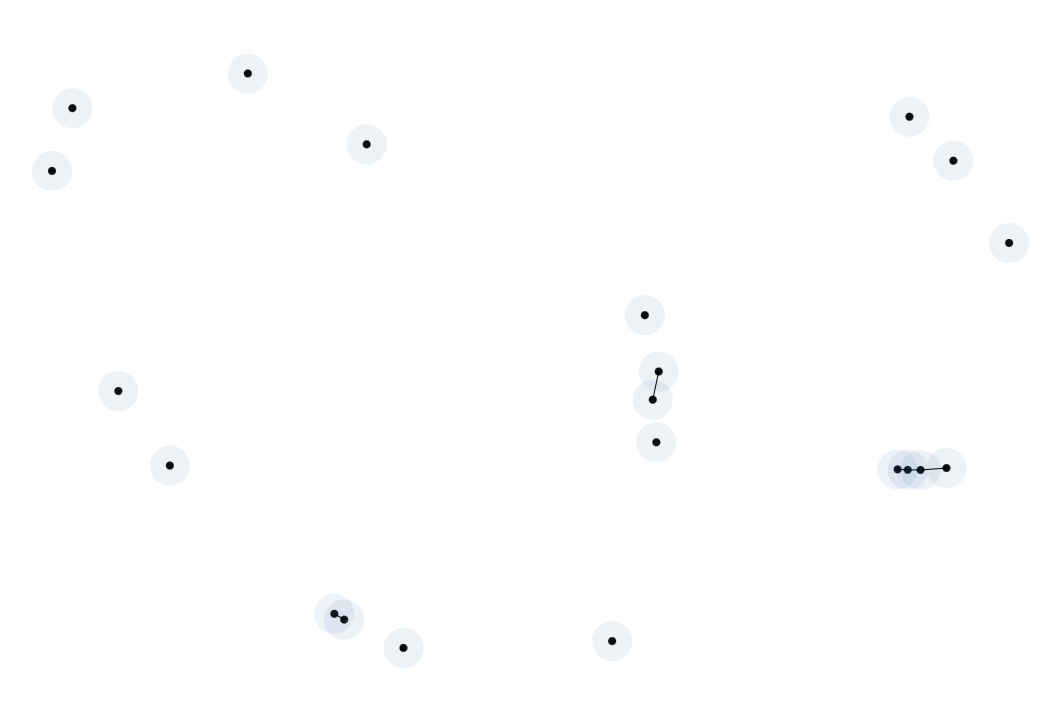}};
     			\node at (4,0) {\includegraphics[width=80pt]{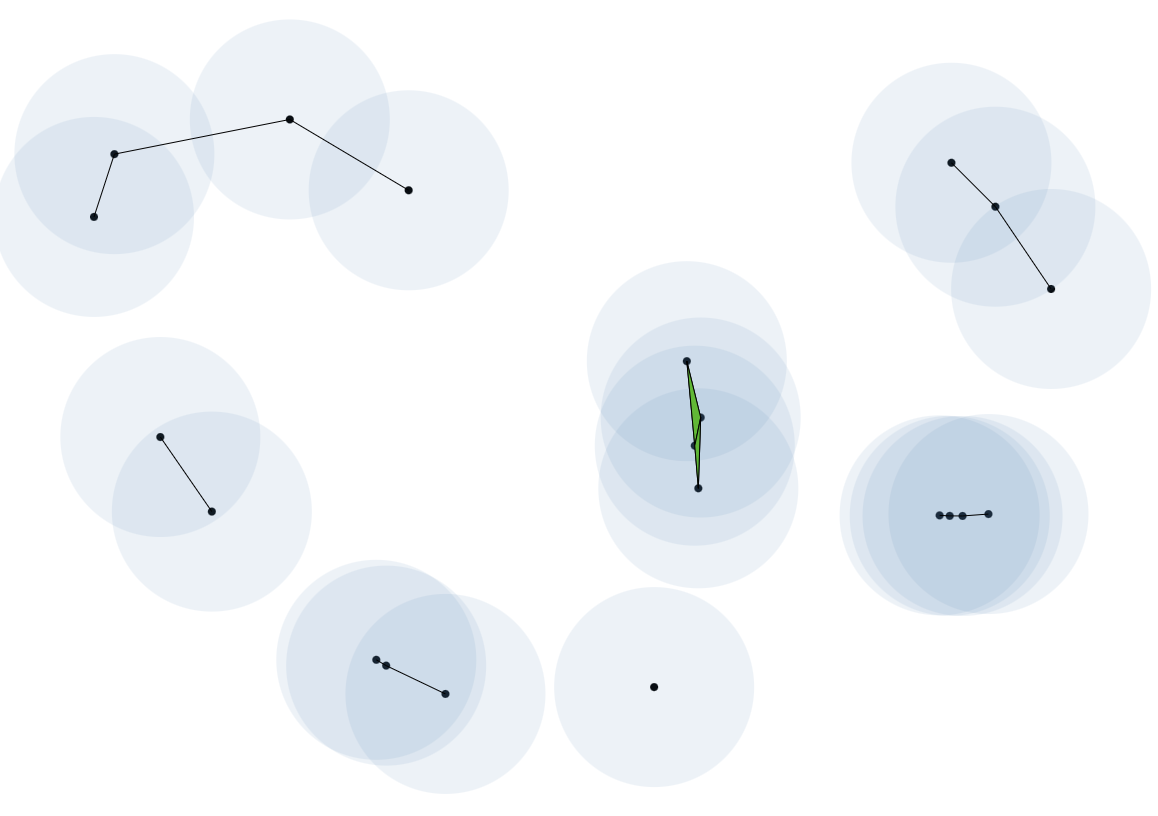}};
     			\node at (8,0) {\includegraphics[width=80pt]{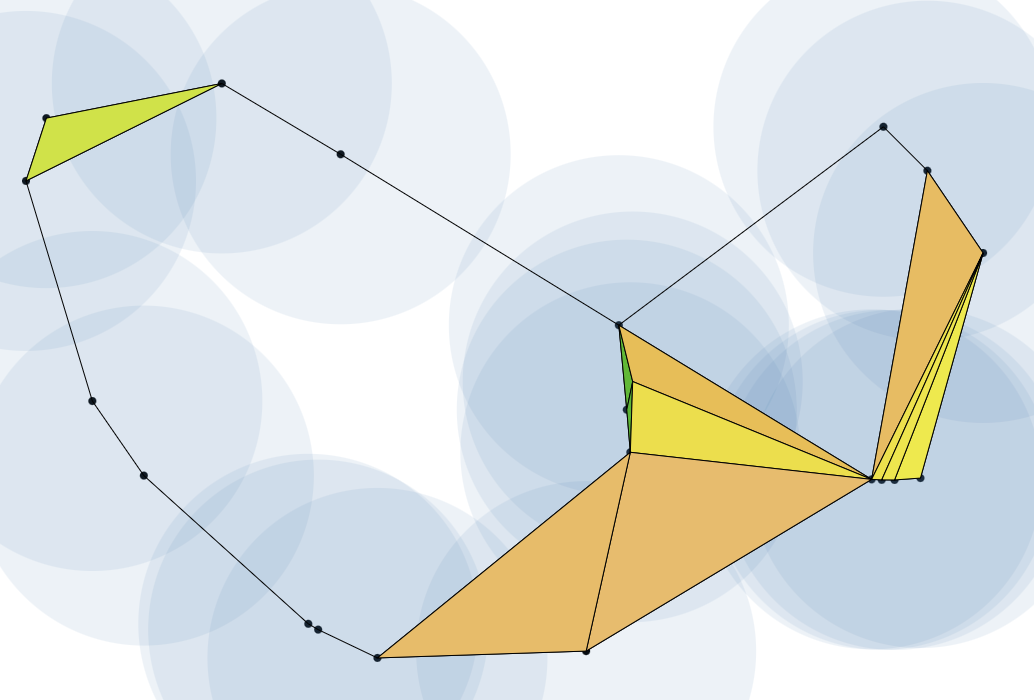}};
     			\node at (12,0) {\includegraphics[width=80pt]{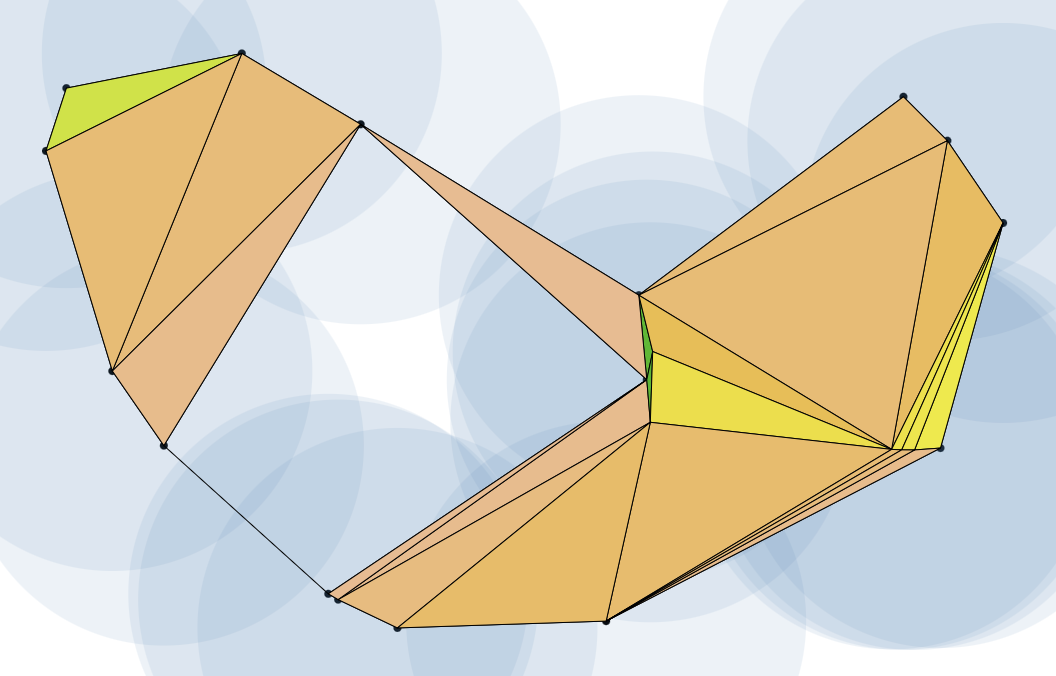}};
     			\node at (6,-4) {\includegraphics[width=200pt]{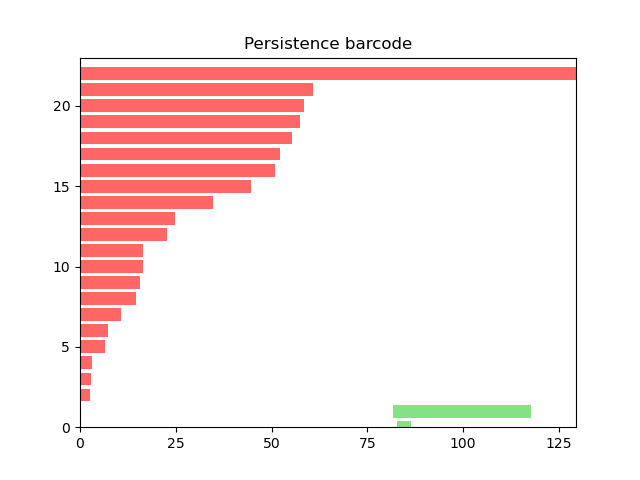}};
     		\end{tikzpicture}
     		\caption{VR complexes and persistent barcodes}
     	\end{figure}
     \end{center}

	\section{Persistent hypergraphs homology}
	
	In this section, we establish the mathematical background for persistent theory of hypergraphs. More details about topology can be found in classical textbooks$^{27,28}$. The base ring $R$ could be any field. For simplicity, we will assume $R$ is of characteristic 2 in the computations.
	
	\subsection{Hypergraphs}
	
	\begin{Definition}
		A \textbf{hypergraph} $\mathcal{H}$ is a pair $(V,E)$ where $V=\{1,2,\dots,|V|\}$ is a finite set and $E\subseteq 2^V$. Elements in $E$ are called \textbf{hyperedges}. Let $\mathcal{H}_n=R$-span$\{e\in E:|e|=n+1\}$ and let $C$ be the simplicial complex $(V,2^V)$, $C_n$ be the collection of all $n$-dim simplices. we view $\mathcal{H}_n$ as a subspace of $C_n$. 
	\end{Definition}
	
	Throughout this subsection, let $\partial_n:C_n\to C_{n-1}$ be the boundary map of $C_*$. Consider the following complexes based on the hypergraph $\mathcal{H}$.
	
	\begin{Definition}
		~
		\begin{enumerate}
			\item The \textbf{infimum chain complex} of $\mathcal{H}$, denoted by $\mathcal{H}_*^{\inf}$, is defined as $\mathcal{H}_n^{\inf}=\partial_n^{-1}(\mathcal{H}_{n-1})\cap \mathcal{H}_n$.
			\item The \textbf{supremum chain complex} of $\mathcal{H}$, denoted by $\mathcal{H}_*^{\sup}$, is defined as $\mathcal{H}_n^{\sup}=\partial_{n+1}(\mathcal{H}_{n+1}) + \mathcal{H}_n$.
			\item The \textbf{associated simplicial complex} $\Delta\mathcal{H}=(V,\Delta_\mathcal{H}(E))$ where $\Delta_\mathcal{H}(E)=\{\sigma\in 2^V:\sigma\subseteq e\text{ for some }e\in E\}$.
			\item The \textbf{lower-associated simplicial complex} $\delta\mathcal{H}=(V,\delta_\mathcal{H}(E))$ where $\delta_\mathcal{H}(E)=\{\sigma\in E:e\in E\text{ for all }e\subseteq\sigma\}$.
		\end{enumerate}
	\end{Definition}

	Let $\Delta \mathcal{H}_*$ and $\delta \mathcal{H}_*$ denote the chain complexes associated to $\Delta \mathcal{H}$ and $\delta \mathcal{H}$ respectively. The boundary maps for $\mathcal{H}^{\inf}_*$, $\mathcal{H}^{\sup}_*$, $\Delta \mathcal{H}_*$ and $\delta \mathcal{H}_*$ are the restricted maps obtained from $\partial_n$. In particular, we denote $\partial_n|_{\mathcal{H}_n^{\inf}}$ and $\partial_n|_{\mathcal{H}_n^{\sup}}$ by $\partial_n^{\inf}$ and $\partial_n^{\sup}$ respectively. The $n$-th homology of $\mathcal{H}^{\inf}_*$, $\mathcal{H}^{\sup}_*$, $\Delta \mathcal{H}_*$ and $\delta \mathcal{H}_*$ are denoted by $H_n^{\inf}$, $H_n^{\sup}$, $\Delta H_n$ and $\delta H_n$ respectively.
	
	Wu and others prove the following insightful result.
	
	\begin{theorem}[Wu and others$^{17}$] For a hypergraph $\mathcal{H}$, one has 
		 $\partial_n(\mathcal{H}_n^{\inf})\subseteq \mathcal{H}_{n-1}^{\inf}$ and $\partial_n(\mathcal{H}_n^{\sup})\subseteq \mathcal{H}_{n-1}^{\sup}$. Hence, $\mathcal{H}_*^{\inf}$ and $\mathcal{H}_*^{\sup}$ are indeed chain complexes. Moreover, $H_n^{\inf}\cong H_n^{\sup}$.
	\end{theorem}

	From their result, $H_n^{\inf}$ is called the \textbf{embedded homology} of $\mathcal{H}$.
	
	It is not hard to see that $\ker(\partial_n^{\inf})=\ker(\partial_n|_{\mathcal{H}_n})$ and $\img(\partial_n^{\sup})=\partial_n(\mathcal{H}_n)$. To better understand the embedded homology, we study $\ker(\partial_n^{\sup})$, $\img(\partial_n^{\inf})$ and their quotient. This quotient space, denoted by $\hat{H}_n=\ker(\partial_n^{\sup})/\img(\partial_{n+1}^{\inf})$ is called the $\hat{H}$ homology of $\mathcal{H}$.
	
	Consider the following commutative diagram
	\begin{center}
		\begin{tikzpicture}
			\node (s) at (3,0) {$0$};
			\node (s0) at (0,0) {$\mathcal{H}^{\sup}_0$};
			\node (s1) at (-3,0) {$\mathcal{H}^{\sup}_1$};
			\node (s2) at (-6,0) {$\mathcal{H}^{\sup}_2$};
			\node (ss) at (-9,0) {$\cdots$};
			\node (i) at (3,-2) {$0$};
			\node (i0) at (0,-2) {$\mathcal{H}^{\inf}_0$};
			\node (i1) at (-3,-2) {$\mathcal{H}^{\inf}_1$};
			\node (i2) at (-6,-2) {$\mathcal{H}^{\inf}_2$};
			\node (ii) at (-9,-2) {$\cdots$};
			
			\draw[black,->] (s0) -- (s);
			\draw[black,->] (s1) -- (s0);
			\draw[black,->] (s2) -- (s1);
			\draw[black,->] (ss) -- (s2);
			\draw[black,->] (i0) -- (i);
			\draw[black,->] (i1) -- (i0);
			\draw[black,->] (i2) -- (i1);
			\draw[black,->] (ii) -- (i2);
			
			\draw[black,right hook->] (i) -- (s);
			\draw[black,right hook->] (i0) -- (s0);
			\draw[black,right hook->] (i1) -- (s1);
			\draw[black,right hook->] (i2) -- (s2);
			
			\node at (-1.5,0.25) {$\partial_1^{\sup}$};
			\node at (-4.5,0.25) {$\partial_2^{\sup}$};
			\node at (-7.5,0.25) {$\partial_3^{\sup}$};
			\node at (-1.5,-1.75) {$\partial_1^{\inf}$};
			\node at (-4.5,-1.75) {$\partial_2^{\inf}$};
			\node at (-7.5,-1.75) {$\partial_3^{\inf}$};
			\node at (-0.25,-1) {$\iota_0$};
			\node at (-3.25,-1) {$\iota_1$};
			\node at (-6.25,-1) {$\iota_2$};
		\end{tikzpicture}
	\end{center}
	where $\iota_n:\mathcal{H}_n^{\inf}\to \mathcal{H}_n^{\sup}$ is the natural embedding.
	
	For each $n$, since $\partial_{n}^{\sup}\circ(\iota_{n}\circ\partial_{n+1}^{\inf})=\iota_{n-1}\circ\partial_{n}^{\inf}\circ\partial_{n+1}^{\inf}=0$, we have the following chain complex
	$$\mathcal{H}^n_*:\cdots \xrightarrow{\partial_{n+3}^{\inf}}\mathcal{H}_{n+2}^{\inf}\xrightarrow{\partial_{n+2}^{\inf}} \mathcal{H}^{\inf}_{n+1} \xrightarrow{\iota_{n}\circ\partial_{n+1}^{\inf}} \mathcal{H}^{\sup}_{n}\xrightarrow{\partial_{n}^{\sup}} \mathcal{H}^{\sup}_{n-1}\xrightarrow{\partial_{n-1}^{\sup}} \cdots.$$
	Under this setting, the $\hat{H}_n$ homology is just the $n$-th homology of $\mathcal{H}^n_*$. In general, $H_n^{\inf}$, $\hat{H}_n$, $\Delta H_n$ and $\delta H_n$ could be mutually distinct.
	
	\subsection{Morphisms and filtrations}
	
	\begin{Definition}
		A \textbf{hypergraph morphism} $\phi:(V_1,E_1)\to (V_2,E_2)$ is a map $\phi:V_1\to V_2$ such that for any hyperedge $e\in E_1$, its image $\{\phi(v):v\in e\}$ is a hyperedge in $E_2$. We say $\phi$ is \textbf{injective} if the map on vertices $\phi:V_1\to V_2$ is injective. In particular, $|e|=|\phi(e)|$ for all hyperedges $e\in E_1$. We say $\phi$ is an \textbf{embedding} if $\phi$ is injective and $\phi(e)=\phi(e')$ implies $e=e'$ for all hyperedges $e,e'\in E_1$.
	\end{Definition}
	
	An injective hypergraph morphism $\phi:\mathcal{H}\to \mathcal{K}$ between hypergraphs $\mathcal{H}$ and $\mathcal{K}$ induces morphisms on chain complexes
	\begin{itemize}
		\item $\phi_*^{\inf}:\mathcal{H}_*^{\inf}\to \mathcal{K}_*^{\inf}$,
		\item $\phi_*^{\sup}:\mathcal{H}_*^{\sup}\to \mathcal{K}_*^{\sup}$,
		\item $\Delta\phi_*:\Delta \mathcal{H}_*\to \Delta \mathcal{K}_*$,
		\item $\delta\phi_*:\delta \mathcal{H}_*\to \delta \mathcal{K}_*$.
	\end{itemize}
	Consequently, $\phi$ induces morphisms of the homology spaces $\Delta H_*$, $\delta H_*$, $H_*^{\inf}$ and $\hat{H}_*$.
	
	\begin{Definition}
		Fix a finite set $V=\{1,2,\dots,|V|\}$. A \textbf{hypergraph filtration} is a map $f:2^V\to\mathbb{N}\cup\{\infty\}$ such that $f(\emptyset)=0$. A hypergraph filtration is a \textbf{simplicial complex filtration} if $f(e)\leq f(e')$ whenever $e\subseteq e'$.
	\end{Definition}
	
	A filtration induces a family of hypergraphs $\{\mathcal{H}^t=(V,E^t):t\in\mathbb{N}\}$ with $E^t=\{e\in 2^V:f(e)\leq t\}$. For $t<r$, one has the natural embedding of hypergraphs $\phi^{t,r}:\mathcal{H}^t\to \mathcal{H}^r$ that further induces morphisms on chain complexes. The $t,r$-\textbf{persistent homology} is defined as
	$$\overline{H}_n^{t,r}=\ker\left(\partial_n|_{\overline{\mathcal{H}}_n^t}\right)/\left(\ker\left(\partial_n|_{\overline{\mathcal{H}}_n^t}\cap\partial_{n+1}\left(\overline{\mathcal{H}}_{n+1}^r\right)\right)\right)$$
	for $\overline{\mathcal{H}}=\Delta \mathcal{H}, \delta \mathcal{H}$ or $\mathcal{H}^{\inf}$. In the other words, let $\left(\overline{\phi}^{t,r}\right)_n^*:\overline{H}_n^t\to\overline{H}_n^r$ be the induced morphism on homology $\overline{H}=\Delta H,\delta H,H^{\inf}$ or $\hat{H}$, then $\overline{H}_n^{t,r}=\img\left(\left(\overline{\phi}^{t,r}\right)_n^*\right)$.
	
	When $f$ is a simplicial complex filtration, $f$ induces the classical persistent homology $H_*^{t,r}$ of simplicial complex as introduced by Zomorodian and Carlsson$^{30}$.
	
	\subsection{Persistence module}
	
	A \textbf{persistence module} is a family of $R$-modules $\{M_i:i\in\mathbb{N}\}$ with $R$-linear maps $\phi_i:M_i\to M_{i+1}$. In this paper, we only consider the case where $R$ is a field of characteristic 2 and $M_i$ are generated by finitely many hyperedges. In this setting, $(M_i,\phi_i)$ forms a persistence module of finite type and has an elegant decomposition.
	
	Each persistence module $(M_i,\phi_i)$ can be identified as a graded $R[t]$-module $\displaystyle\bigoplus_{i\geq 0}M_i$ where the action of $t$ is given by
	$$t\cdot(m_0,m_1,m_2,m_3,\dots)=(0,\phi_0(m_0),\phi_1(m_1),\phi_2(m_2),\dots).$$
	
	Since $R[t]$ is a principal ideal domain, we have the following decomposition, followed from the classical structural theorem in commutative algebra.
	
	\begin{theorem}
		$$\bigoplus_{i\geq 0}M_i=\left(\bigoplus_{i=1}^n \Sigma^{\alpha_i} R[t] \right)\oplus \left(\bigoplus_{j=1}^m \Sigma^{\beta_j} R[t]/(t^{\gamma_j}) \right)$$
		where $\Sigma^a$ means shifting the grading up by $a$.
	\end{theorem}

	In the decomposition, each summand can be interpreted by a pair $(i,j)\in \mathbb{N}\times\mathbb{N}\cup\{\infty\}$ as $\Sigma^{\alpha_i} R[t]\mapsto (\alpha_i,\infty)$ and $\Sigma^{\beta_j} R[t]/(t^{\gamma_j})\mapsto (\beta_j,\beta_j+\gamma_j)$. The multiset $\{(\alpha_i,\infty),(\beta_j,\beta_j+\gamma_j):1\leq i\leq n,1\leq j\leq m\}$ is called the \textbf{barcodes} representation of the persistent module $(M_i,\phi_i)$.

     \section{Computational experiments and data analysis}
	
	\subsection{Dataset}
	The data we use in this paper is provided by SocioPatterns$^{30}$ sensing program. As is shown in Table 1, the data consists of face-to-face interactions of the following social groups
	\begin{enumerate}
		\item 13 Guinea baboons living at CNRS primate center$^{25}$;
		\item 86 individuals from a village in rural Malawi$^{26}$;
		\item 403 participants of a scientific conference$^{24}$;
		\item 217 people at a workplace$^{24}$;
		\item 327 students from a high school in Marseilles$^{22}$.
	\end{enumerate}
	
	\begin{table}\label{tab:data}
		\begin{tabular}{cccc}
			\hline
			data sets & number of individuals & total number of contacts & duration of data collection \\
			\hline
			Baboon & 13 & 63095 & 4 weeks\\
			Malawi & 86 & 102293 & 2 weeks\\
			Conference & 403 & 70261 & 2 days\\
			Workplace & 217 & 78249 & 2 weeks\\
			High school & 327 & 188508 & 1 week\\
			\hline
			& & &
		\end{tabular}
		\caption{Characteristics of the data sets.}
	\end{table}
	
	All of these face-to-face contact data sets are collected and recorded in unsupervised fashion by using wearable sensors called Radio-Frequency Identification Devices (RFID). The sensors send signal once every 20 seconds, and they can receive signals from other devices within the radius of 1.5 meter. The sensors are weared on the chest of each individuals so that their devices will receive and record signals when they are facing each other. If a pair of devices both receive signal from each other within a period of 20 seconds, we say the two individuals are in face-to-face contact in that period. It is shown that wearable sensors detect face-to-face interactions with more than 99\% accuracy.

	\subsection{Persistent hypergraph model}
	For each set of data, we construct a persistent hypergraph. For each hyperedge $e$ consisting of $k$ individuals, we record $T_e$, the total number of periods during which every pair among them are in face-to-face contact. In the other words, $T_e$ is the total length of time that these $k$ people are in a group meeting. It should be noted that if there is a group meeting, only the hyperedge consisting of the whole group counts. We ignore all the subgroups to avoid repeated counting.

	The filtration for a hyperedge $e$ is given by $M - \log(T_e)$ where $M=\log(\max\{T_{f}: f\text{ is hyperedge}\})$ is the fixed constant ensuring the non-negativity of the filtration. By convention, the filtration of $e$ is set to infinity if $T_e=0$.
	
	To better explain our model, we give an easy but complete example. Suppose our data consists of six individuals. The interactions are given in Figure 2 in terms of graphs i.e. in $t_k$, an edge $(i,j)$ means there is an interaction between $i$ and $j$ at time $k$.
	
	\begin{center}
		\begin{figure}\label{fig:fakedata}
			\begin{tikzpicture}
				\node at (0,0) {\scalebox{0.7}{\begin{tikzpicture}
						\node (A) at (-1,1.464) {A};
						\node (B) at (1,1.464) {B};
						\node (C) at (2,0) {C};
						\node (D) at (1,-1.464) {D};
						\node (E) at (-1,-1.464) {E};
						\node (F) at (-2,0) {F};
						
						\draw[black,thick] (A) -- (B) -- (C) -- (A);
						\draw[black,thick] (D) -- (E) -- (F) -- (D);
						
						\node at (0,-2.5) {\large $t_1$};
				\end{tikzpicture}}};
				\node at (4,0) {\scalebox{0.7}{\begin{tikzpicture}
							\node (A) at (-1,1.464) {A};
							\node (B) at (1,1.464) {B};
							\node (C) at (2,0) {C};
							\node (D) at (1,-1.464) {D};
							\node (E) at (-1,-1.464) {E};
							\node (F) at (-2,0) {F};
							
							\draw[black,thick] (A) -- (B) -- (C) -- (A);
							\draw[black,thick] (F) -- (D);
							
							\node at (0,-2.5) {\large $t_2$};
				\end{tikzpicture}}};
				\node at (8,0) {\scalebox{0.7}{\begin{tikzpicture}
							\node (A) at (-1,1.464) {A};
							\node (B) at (1,1.464) {B};
							\node (C) at (2,0) {C};
							\node (D) at (1,-1.464) {D};
							\node (E) at (-1,-1.464) {E};
							\node (F) at (-2,0) {F};
							
							\draw[black,thick] (A) -- (B);
							\draw[black,thick] (D) -- (C);
							
							\node at (0,-2.5) {\large $t_3$};
				\end{tikzpicture}}};
				\node at (0,-4) {\scalebox{0.7}{\begin{tikzpicture}
							\node (A) at (-1,1.464) {A};
							\node (B) at (1,1.464) {B};
							\node (C) at (2,0) {C};
							\node (D) at (1,-1.464) {D};
							\node (E) at (-1,-1.464) {E};
							\node (F) at (-2,0) {F};
							
							\draw[black,thick] (C) -- (A);
							\draw[black,thick] (F) -- (D);
							
							\node at (0,-2.5) {\large $t_4$};
				\end{tikzpicture}}};
				\node at (4,-4) {\scalebox{0.7}{\begin{tikzpicture}
							\node (A) at (-1,1.464) {A};
							\node (B) at (1,1.464) {B};
							\node (C) at (2,0) {C};
							\node (D) at (1,-1.464) {D};
							\node (E) at (-1,-1.464) {E};
							\node (F) at (-2,0) {F};
							
							\draw[black,thick] (B) -- (C);
							\draw[black,thick] (D) -- (E) -- (F) -- (D);
							
							\node at (0,-2.5) {\large $t_5$};
				\end{tikzpicture}}};
				\node at (8,-4) {\scalebox{0.7}{\begin{tikzpicture}
							\node (A) at (-1,1.464) {A};
							\node (B) at (1,1.464) {B};
							\node (C) at (2,0) {C};
							\node (D) at (1,-1.464) {D};
							\node (E) at (-1,-1.464) {E};
							\node (F) at (-2,0) {F};
							
							\draw[black,thick] (A) -- (B) -- (C) -- (A);
							\draw[black,thick] (F) -- (D);
							
							\node at (0,-2.5) {\large $t_6$};
				\end{tikzpicture}}};
				\node at (0,-8) {\scalebox{0.7}{\begin{tikzpicture}
							\node (A) at (-1,1.464) {A};
							\node (B) at (1,1.464) {B};
							\node (C) at (2,0) {C};
							\node (D) at (1,-1.464) {D};
							\node (E) at (-1,-1.464) {E};
							\node (F) at (-2,0) {F};
							
							\draw[black,thick] (A) -- (B);
							\draw[black,thick] (D) -- (C);
							
							\node at (0,-2.5) {\large $t_7$};
				\end{tikzpicture}}};
				\node at (4,-8) {\scalebox{0.7}{\begin{tikzpicture}
							\node (A) at (-1,1.464) {A};
							\node (B) at (1,1.464) {B};
							\node (C) at (2,0) {C};
							\node (D) at (1,-1.464) {D};
							\node (E) at (-1,-1.464) {E};
							\node (F) at (-2,0) {F};
							
							\draw[black,thick] (A) -- (B);
							\draw[black,thick] (F) -- (C);
							
							\node at (0,-2.5) {\large $t_8$};
				\end{tikzpicture}}};
				\node at (8,-8) {\scalebox{0.7}{\begin{tikzpicture}
							\node (A) at (-1,1.464) {A};
							\node (B) at (1,1.464) {B};
							\node (C) at (2,0) {C};
							\node (D) at (1,-1.464) {D};
							\node (E) at (-1,-1.464) {E};
							\node (F) at (-2,0) {F};
							
							\draw[black,thick] (A) -- (B);
							\draw[black,thick] (F) -- (C);
							
							\node at (0,-2.5) {\large $t_9$};
				\end{tikzpicture}}};
			\end{tikzpicture}
		\caption{Interactions between six individuals, an edge $(i,j)$ in $t_k$ means an interaction takes place between $i$ and $j$ at time $k$.}
		\end{figure}
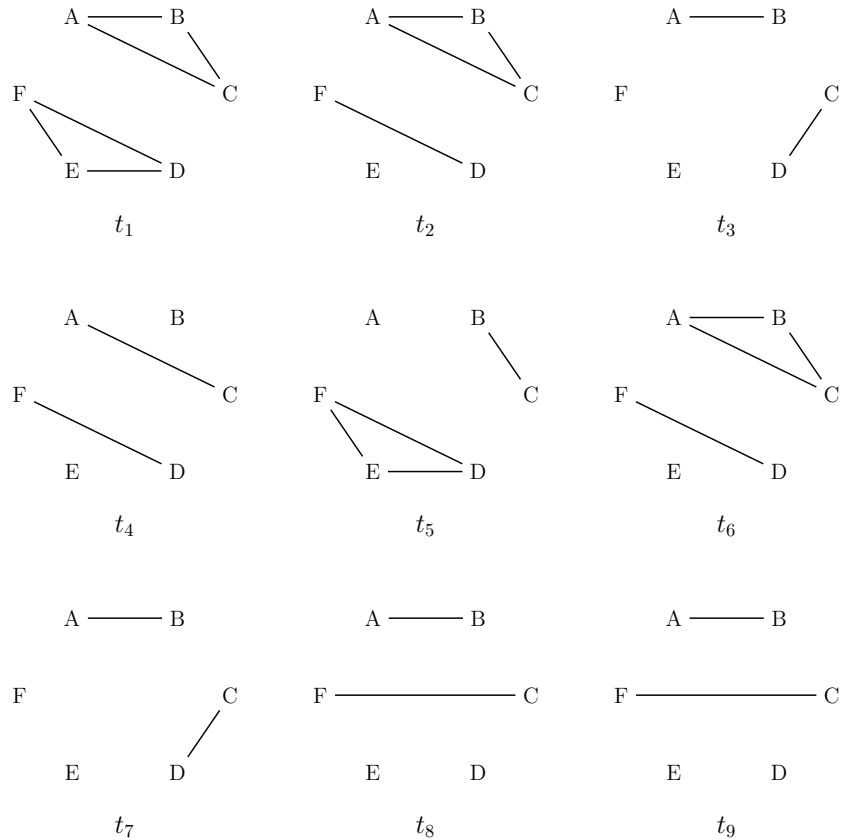
	\end{center}
	
	Then, the filtration function and persistent hypergraphs are given in Figure 3.
	
	\begin{center}
		\begin{figure}\label{fig:hypergraph}
			\begin{tabular}{|c|c|c|c|c|c|c|c|c|}
				\hline
				$e$ & AB & DF & CD & CF & AC & BC & ABC & DEF \\
				\hline
				$T_e$ & 4 & 3 & 2 & 2 & 1 & 1 & 3 & 2 \\
				\hline
				$f(e)$ & 0 & $\log 4/3$ & $\log 2$ & $\log 2$ & $\log 4$ & $\log 4$ & $\log 4/3$ & $\log 2$ \\
				\hline
			\end{tabular}
		
			\begin{tikzpicture}
				\node at (0,0) {\includegraphics[width=80pt]{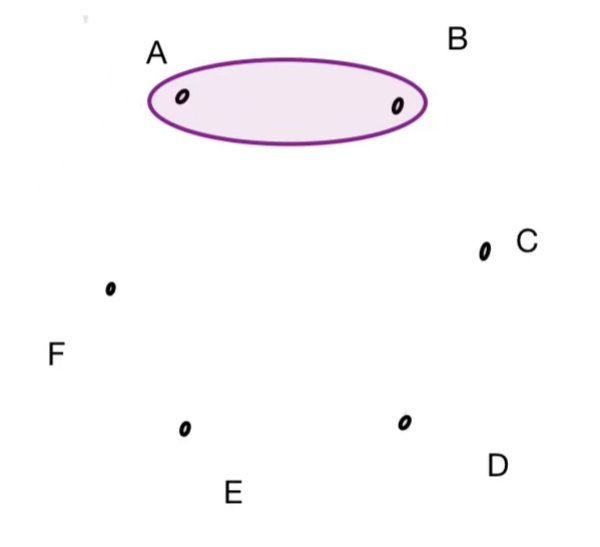}};
				\node at (4,0) {\includegraphics[width=80pt]{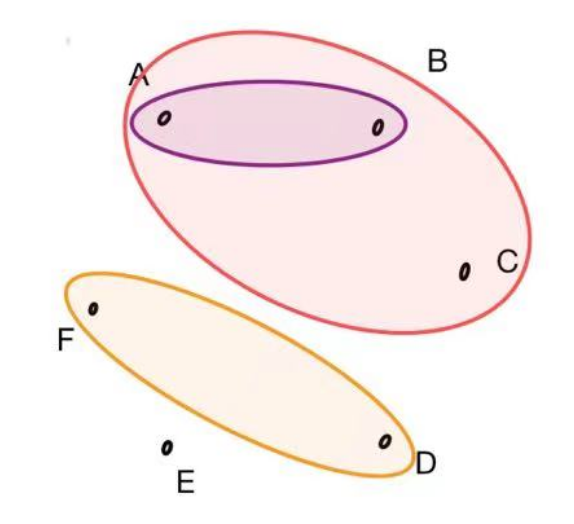}};
				\node at (8,0) {\includegraphics[width=80pt]{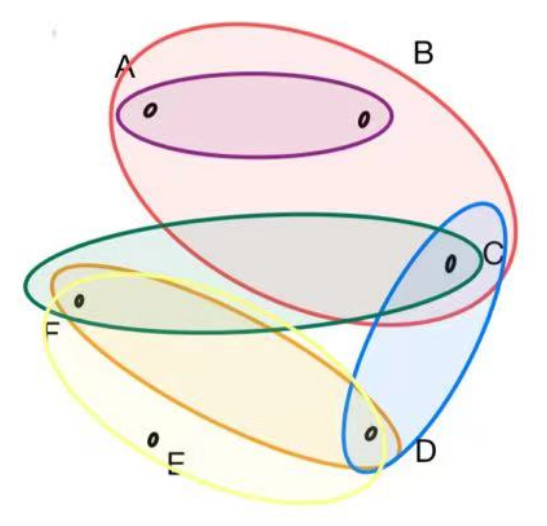}};
				\node at (12,0) {\includegraphics[width=80pt]{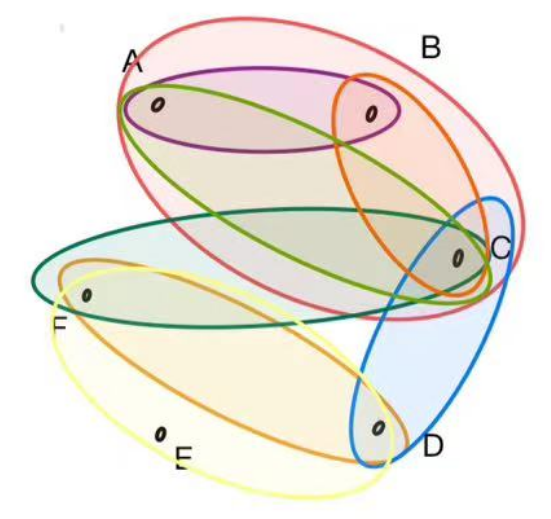}};
				
				\node at (0,-1.7) {$t=0$};
				\node at (4,-1.7) {$t=\log 4/3$};
				\node at (8,-1.7) {$t=\log 2$};
				\node at (12,-1.7) {$t=\log 4$};
				
				\node at (3.6,-2.5) {Barcodes for $H_0^{\inf}$ in dark red, and $H_1^{\inf}$ in dark blue};
				\filldraw[dr] (0,-2.9) -- (13,-2.9) -- (13,-3.0) -- (0,-3.0) -- (0,-2.9);
				\filldraw[dr] (13,-2.8) -- (13,-3.1) -- (13.3,-2.95);
				\filldraw[dr] (0,-3.1) -- (12,-3.1) -- (12,-3.2) -- (0,-3.2) -- (0,-3.1);
				\filldraw[dr] (0,-3.3) -- (12,-3.3) -- (12,-3.4) -- (0,-3.4) -- (0,-3.3);
				\filldraw[dr] (0,-3.5) -- (8,-3.5) -- (8,-3.6) -- (0,-3.6) -- (0,-3.5);
				\filldraw[dr] (0,-3.7) -- (8,-3.7) -- (8,-3.8) -- (0,-3.8) -- (0,-3.7);
				\filldraw[db] (8,-3.9) -- (13,-3.9) -- (13,-4) -- (8,-4) -- (8,-3.9);
				\filldraw[db] (13,-3.8) -- (13,-4.1) -- (13.3,-3.95);
				
				\node at (2.5,-4.4) {Additional barcodes for $\hat{H}_1$ in dark blue};
				\filldraw[db] (4,-4.8) -- (12,-4.8) -- (12,-4.9) -- (4,-4.9) -- (4,-4.8);
				\filldraw[db] (8,-5) -- (13,-5) -- (13,-5.1) -- (8,-5.1) -- (8,-5);
				\filldraw[db] (13,-4.9) -- (13,-5.2) -- (13.3,-5.05);
			\end{tikzpicture}
			\captionof{figure}{Hyperedges $e$, duration of interactions $T_e$, filtrations $f(e)$, persistent hypergraphs and persistent barcodes}
		\end{figure}
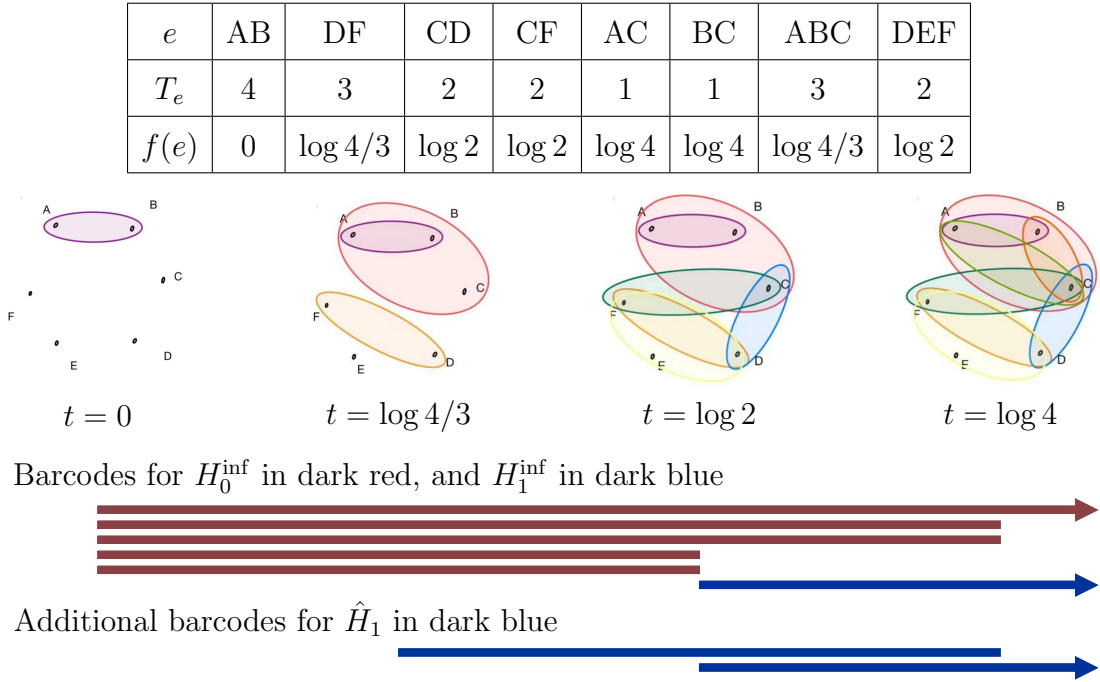
	\end{center}
	
	After these preparations, persistent embedded homology and $\hat{H}$ homology of dimension 0 and 1 are computed using our algorithm.
	
	\subsection{Algorithm}
	
	In this section, we provide an algorithm of computing barcodes for hypergraphs from a given filtration. This algorithm is inspired by the one for simplicial complex introduced by Zomorodian and Carlsson$^{30}$. When the filtration satisfies the condition of a persistent simplicial complex, our algorithm yields the classical barcodes for simplicial complex.
	
	Suppose we fix a filtration $f:2^V\to\mathbb{N}\cup\{\infty\}$, we compute the dimension $k$ persistent $H^{\inf}$ and $\hat{H}$ barcodes simultaneously as follows:
	
	\begin{itemize}
		\item Step 1. Let $E_k=\{e\in 2^V:|e|=k+1\}$ be the set of hyperedges of dimension $k$. Sort $E_k$ decreasingly according to the filtration. Let $E_{k+1}=\{e\in 2^V:|e|=k+2,f(e)<\infty\}$ be the set of hyperedges of dimension $k+1$ whose filtration is less than infinity. Sort $E_{k+1}$ increasingly according to the filtration.
		\item Step 2. Construct a boundary matrix $M=(M_{ij})$ whose rows correspond to $E_k$ and columns correspond to $E_{k+1}$. The entry $M_{ij}=1$ if the $i$-th hyperedge in $E_k$ is contained in the $j$-th hyperedge in $E_{k+1}$, and $M_{ij}=0$ otherwise.
		\item Step 3. Compute pivots from boundary matrix $M$ using the pseudo-codes
		
		\begin{tikzpicture}
			\node at (-0.8,0.5) {\textbf{end}};
			\node at (-0.17,1) {\textbf{end}};
			\node at (0.42,1.5) {\textbf{end}};
			\node at (4.45,2) {\textbf{change} column $k$ to column $k$ $-$ column $j$;};
			\node at (3.6,2.5) {\textbf{if} $k$ is not a pivot column and $M_{i,k}=1$};
			\node at (2.7,3) {\textbf{for} $k=j+1$ \textbf{to} number of columns};
			\node at (1.18,3.5) {\textbf{add} $(i,j)$ to Pivot;};
			\node at (5.5,4) {\textbf{find} the smallest $j$ such that $j$ is not a pivot column and $M_{i,j}=1$;};
			\node at (1.35,4.5) {\textbf{for} $i=1$ \textbf{to} number of rows};
			\node at (0.6,5) {\textbf{initialize} Pivot$=\emptyset$;};
			\node at (0,5.7) {\textbf{ComputePivot}($M$)};
			
			\node at (-1.8,0.5) {10};
			\node at (-1.8,1) {9};
			\node at (-1.8,1.5) {8};
			\node at (-1.8,2) {7};
			\node at (-1.8,2.5) {6};
			\node at (-1.8,3) {5};
			\node at (-1.8,3.5) {4};
			\node at (-1.8,4) {3};
			\node at (-1.8,4.5) {2};
			\node at (-1.8,5) {1};
			
			\draw[thick,black] (-2.2,0.2) -- (11.6,0.2);
			\draw[thick,black] (-2.2,5.3) -- (11.6,5.3);
			\draw[thick,black] (-1.4,0.4) -- (-1.4,5.1);
			
			\node at (0,6.1) {};
		\end{tikzpicture}
		
		\item Step 4. For each row, record an ordered pair (birth, death) where birth is the filtration of the hyperedge corresponding to the row. If there is a pivot at column $j$, then death is the filtration of the hyperedge corresponding to column $j$. If the row has no pivot, then death is set to be infinity.
		\item Step 5. For each ordered pair (birth, death), if birth $<$ death, then the pair corresponds to a barcode in persistent embedded homology $H_n^{\inf}$. If birth $>$ death, then the pair (death, birth) corresponds to a barcode in persistent homology $\hat{H}$.
	\end{itemize}
	The complexity of our algorithm is at the level of $n^{k+2}$ where $n$ is the number of individuals and $k$ is the dimension. 
	This is sufficiently efficient for our purpose. For the network of a conference in which there are 403 individuals and 70261 face-to-face contacts, the computation for dimension 0 and 1 barcodes is done within a few seconds on a regular personal computer.The software used to compute the results in this section is available at
	\url{https://github.com/Gao-Yaru/Interaction-Persistent-Hypergraph}.
	
	There are certainly ways to improve this algorithm. For example, when constructing the boundary matrix, instead of taking all hyperedges for rows, we only need to take those that are subsets of existing column hyperedges. But we will not discuss it further in this paper.
	
	\subsection{Results: 0-dimensional  persistent hypergraph homology}
	
	We begin by analyzing results in dimension 0. In principle, dimension 0 barcodes always begin at 0. Those bars ending at infinity correspond to connected components in the connectivity graph, that is, groups of individuals such that any pair of individuals from different groups never conduct face-to-face interaction. The persistent hypergraph model provides an efficient computation on the number of connected components. Among data sets we use, Malawi has two connected components, meaning that the local residents form two groups, and there is no contact between the two groups during the two weeks of data collection. Our computation also show that all of the other data sets are connected. Figure 4 provides a comparison between barcodes.
	
	\begin{figure}[htbp!]\label{fig:Bar0}
		\begin{center}
			\begin{tikzpicture}
				\node at (0,0) {\includegraphics[width = 80pt]{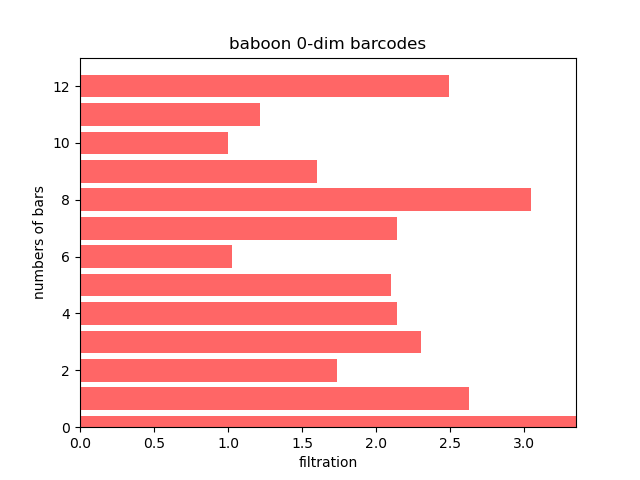}};
				\node at (3,0) {\includegraphics[width = 80pt]{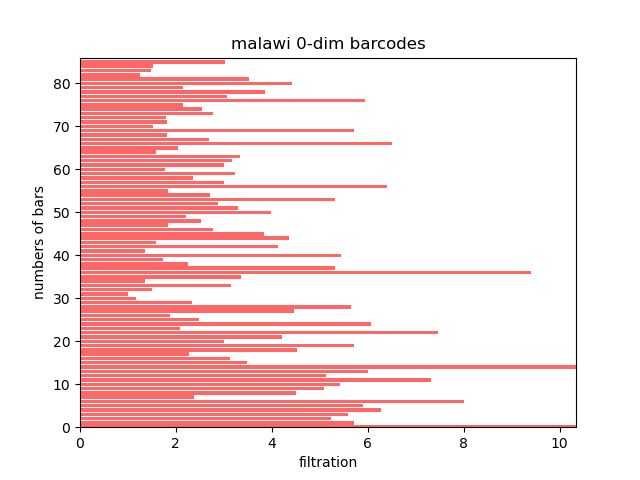}};
				\node at (6,0) {\includegraphics[width = 80pt]{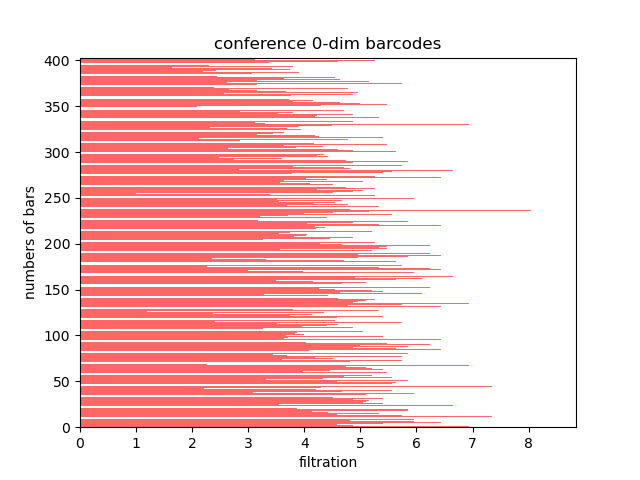}};
				\node at (9,0) {\includegraphics[width = 80pt]{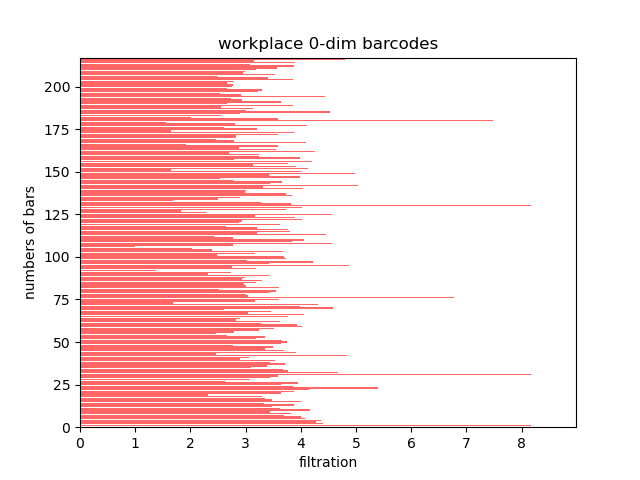}};
				\node at (12,0) {\includegraphics[width = 80pt]{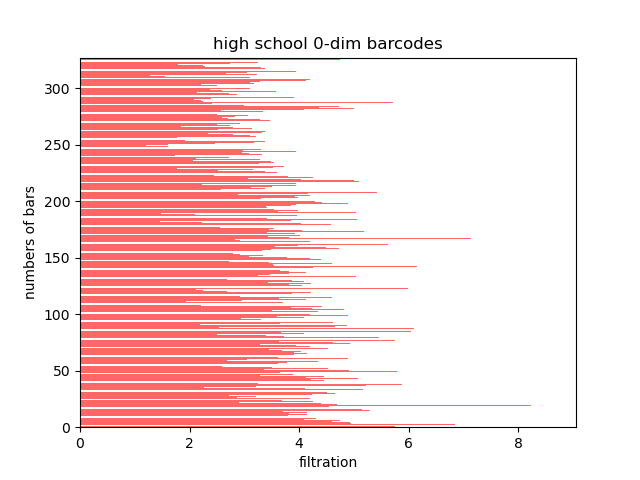}};
				\node at (0,-1.5) {Baboon};
				\node at (3,-1.5) {Malawi};
				\node at (6,-1.5) {Conference};
				\node at (9,-1.5) {Workplace};
				\node at (12,-1.5) {High school};
			\end{tikzpicture}
		\end{center}
		\caption{Dimension 0 barcodes of the datasets.}
	\end{figure}
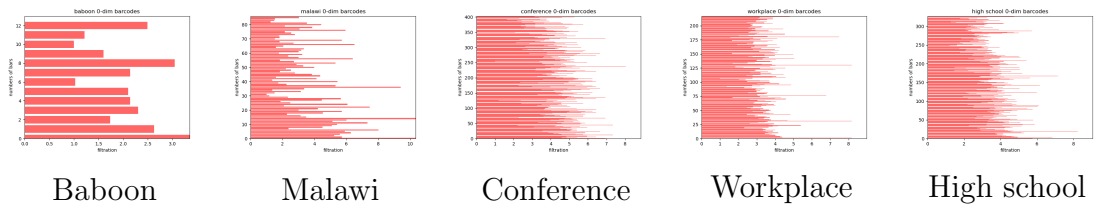
	
	While we can read the number of mutually unrelated groups from persistent barcodes, they do not contain the information of the exact members of the groups. Finding the members requires additional efforts.

		\subsection{Results: 1-dimensional  persistent hypergraph homology}
	
	Much more interesting results appear in dimension 1. In this case, there are two types of barcodes: the persistent $H^{\inf}$ and $\hat{H}$ barcodes. Theoretically, persistent $H^{\inf}$ barcodes correspond to holes in the connectivity graphs, that is a loop whose interior is not filled. In particular, a bar ending at infinity means there exists a group of individuals, say 3 individuals for example, that experience pair-wise face-to-face interaction, but they never meet together simultaneously. On the other hand, an $\hat{H}$ bar corresponds to an anti-hole i.e. a filled interior with some missing edges on its boundary. In particular, an $\hat{H}$ bar ending at infinity means there is a group meeting of 3 individuals among which some pairs never interact privately. As examples, we show these barcodes calculated from the datasets in Figure 5.
	
	\begin{figure}[htbp!]\label{fig:Bar1}
		\begin{center}
			\begin{tikzpicture}
				\node at (0,0) {\includegraphics[width = 80pt]{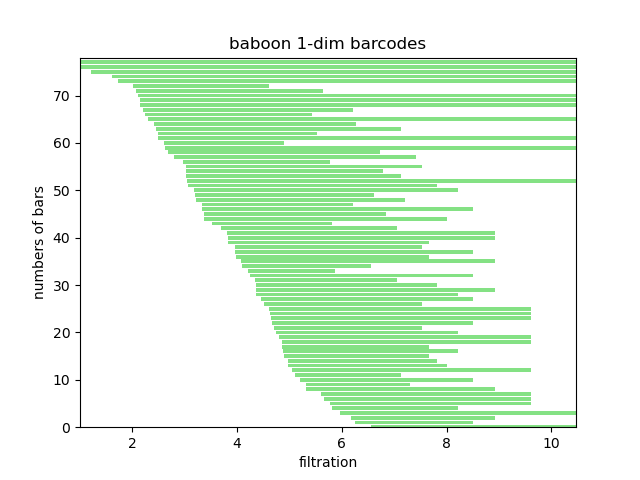}};
				\node at (3,0) {\includegraphics[width = 80pt]{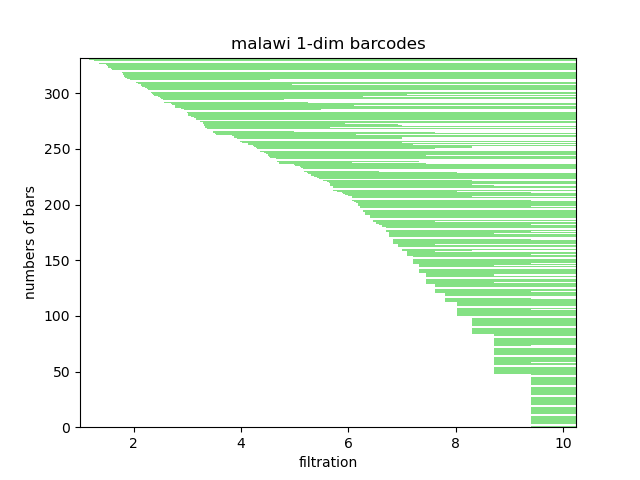}};
				\node at (6,0) {\includegraphics[width = 80pt]{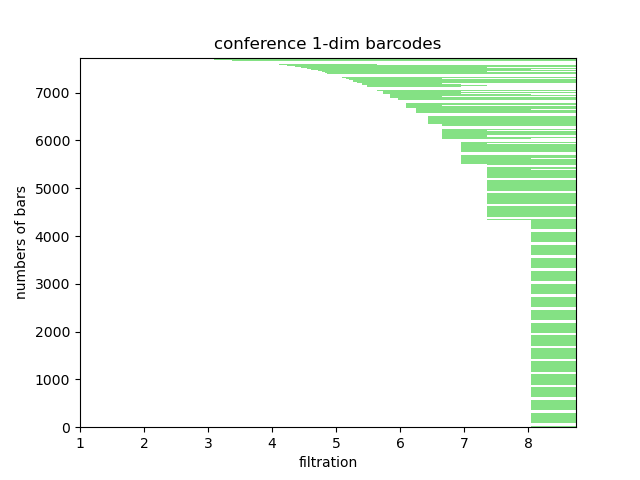}};
				\node at (9,0) {\includegraphics[width = 80pt]{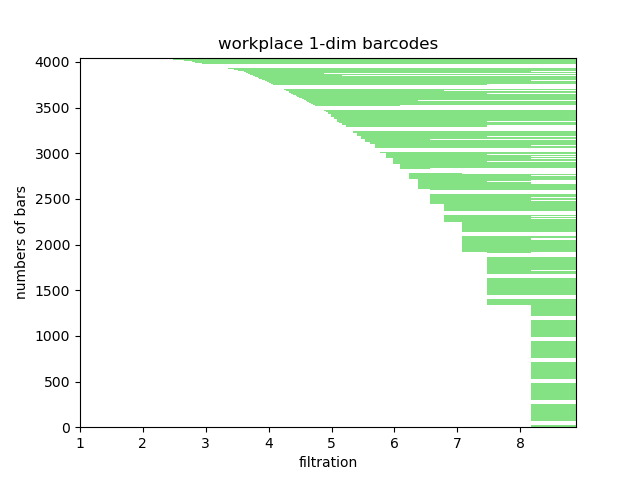}};
				\node at (12,0) {\includegraphics[width = 80pt]{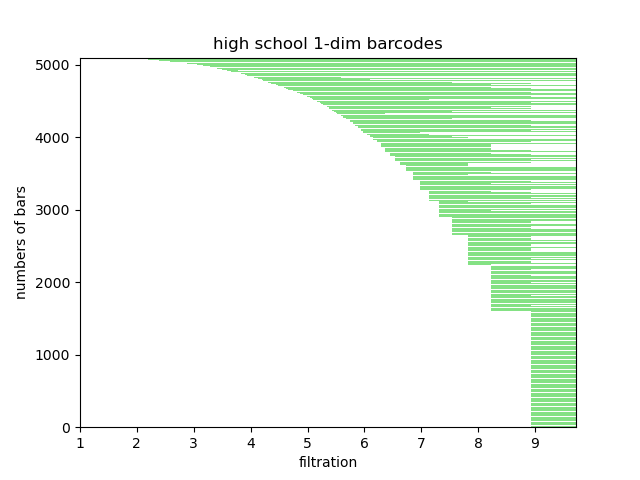}};
				\node at (0,-2.5) {\includegraphics[width = 80pt]{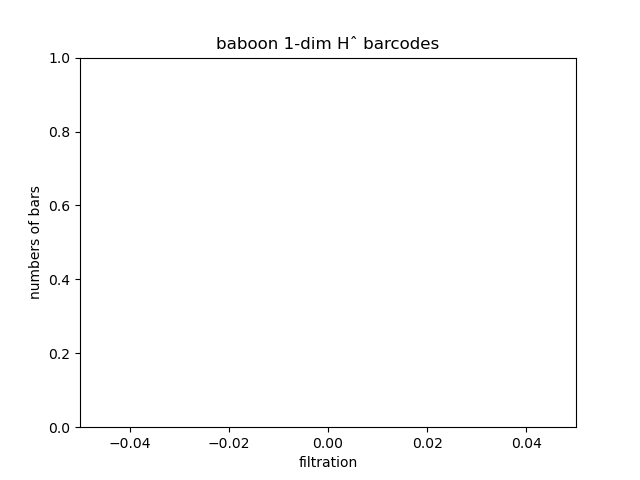}};
				\node at (3,-2.5) {\includegraphics[width = 80pt]{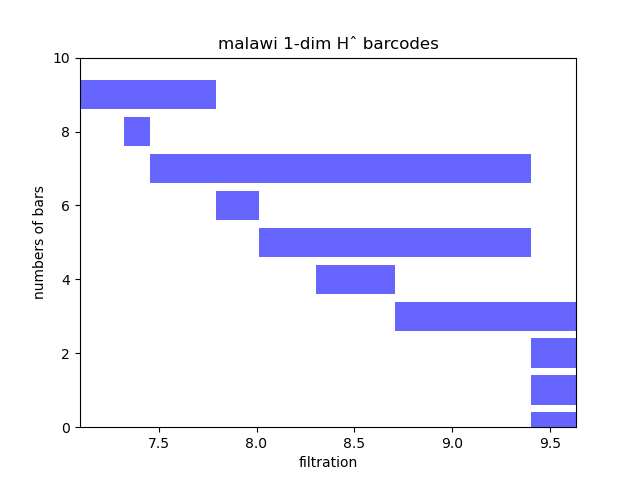}};
				\node at (6,-2.5) {\includegraphics[width = 80pt]{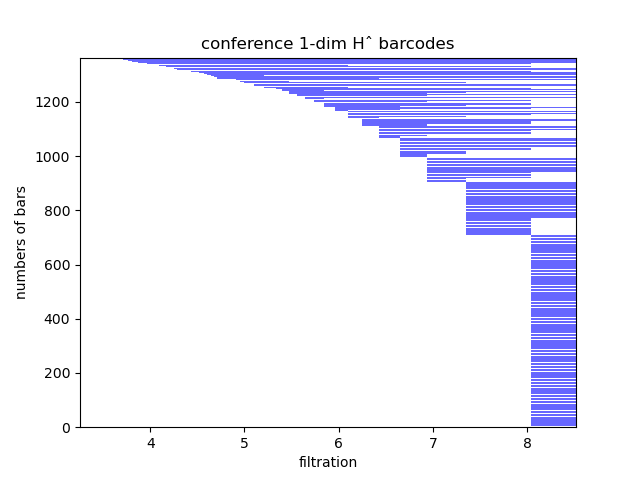}};
				\node at (9,-2.5) {\includegraphics[width = 80pt]{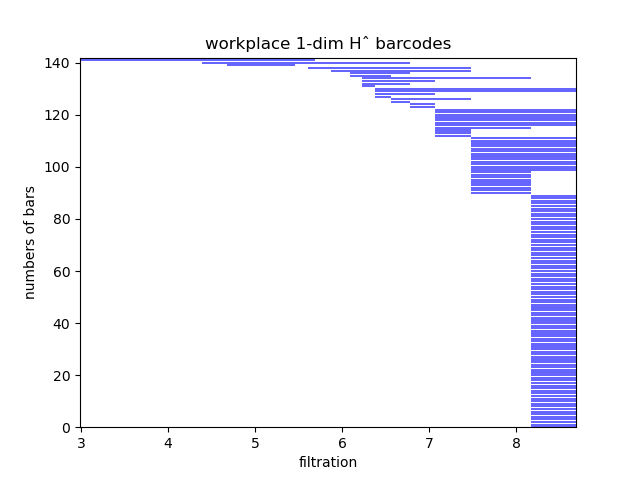}};
				\node at (12,-2.5) {\includegraphics[width = 80pt]{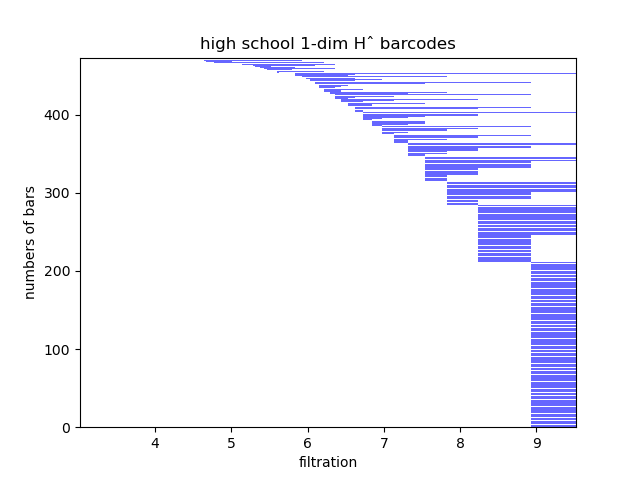}};
				\node at (0,-4) {Baboon};
				\node at (3,-4) {Malawi};
				\node at (6,-4) {Conference};
				\node at (9,-4) {Workplace};
				\node at (12,-4) {High school};
			\end{tikzpicture}
		\end{center}
		\caption{Dimension 1 persistent $H^{\inf}$ (green) and $\hat{H}$ (blue) barcodes of the datasets.}
	\end{figure}
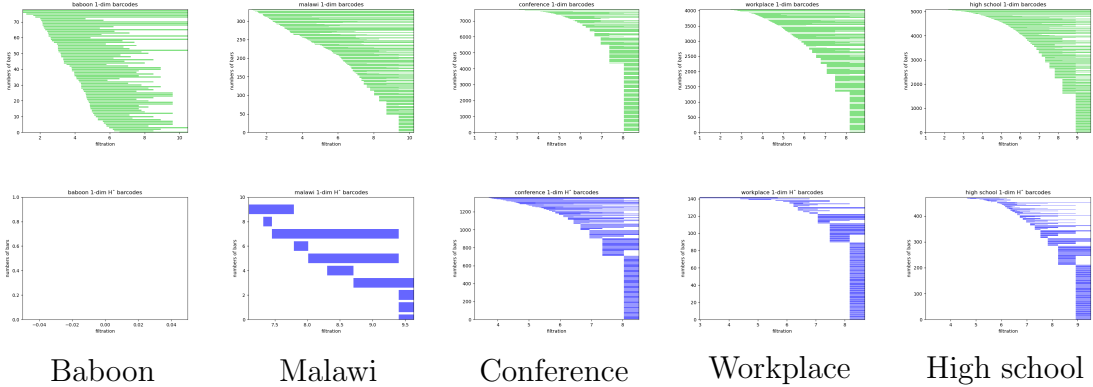
	
	Table 2 provides computation results of the number of persistent $H^{\inf}$ and $\hat{H}$ barcodes and those ending at infinity. We emphasize that the number of barcodes ending at infinity is a new and useful statistic quantity in analysis of population behaviors. More precisely, we focus on the proportions $n/(N+\hat{N})$ and $\hat{n}/(N+\hat{N})$ of the bars ending at infinity among all bars. These results are visualized in Figure 6.
	
	\begin{table}[htbp!]\label{tab:result}
		\begin{tabular}{ccccc}
			\hline
			data sets & $N$ & $n$ & $\hat{N}$ & $\hat{n}$ \\
			\hline
			Baboon & 78 & 14 & 0 & 0\\
			Malawi & 332 & 270 & 10 & 4\\
			Conference & 7733 & 7290 & 1363 & 1072\\
			Workplace & 4048 & 3537 & 142 & 112\\
			High school & 5096 & 3797 & 473 & 281\\
			\hline
			& & & &
		\end{tabular}
		\caption{Statistics on barcodes of dimension 1: $N$ is the number of persistent $H^{\inf}$ barcodes, $n$ is the number of persistent $H^{\inf}$ barcodes ending at infinity, $\hat{N}$ is the number of persistent $\hat{H}$ barcodes, $\hat{n}$ is the number of persistent $\hat{H}$ barcodes ending at infinity.}
	\end{table}
	
	Roughly speaking, a high proportion of persistent barcodes ending at infinity reflects that a high proportion of the population tends to contact privately rather than communicate in groups. In contrast, higher proportion of $\hat{H}$ barcodes ending at infinity represents that population are more open to strangers in the sense that individuals are more active in participating in group meetings with those they have not contacted.
	
	\begin{figure}[htbp!]\label{fig:result}
		\begin{center}
			\scalebox{1}{
				\begin{tikzpicture}
					\filldraw[ForestGreen] (1,0.897) circle (3pt) {};
					\filldraw[ForestGreen] (4,3.947) circle (3pt) {};
					\filldraw[ForestGreen] (7,4.007) circle (3pt) {};
					\filldraw[ForestGreen] (10,4.221) circle (3pt) {};
					\filldraw[ForestGreen] (13,3.409) circle (3pt) {};
					
					\filldraw[dr] (1,0) circle (3pt) {};
					\filldraw[dr] (4,0.146) circle (3pt) {};
					\filldraw[dr] (7,0.589) circle (3pt) {};
					\filldraw[dr] (10,0.134) circle (3pt) {};
					\filldraw[dr] (13,0.252) circle (3pt) {};
					
					\node at (0.8,1.2) {\scalebox{0.8}{$0.897$}};
					\node at (4,4.2) {\scalebox{0.8}{$3.947$}};
					\node at (7,4.3) {\scalebox{0.8}{$4.007$}};
					\node at (10,4.55) {\scalebox{0.8}{$4.221$}};
					\node at (13,3.7) {\scalebox{0.8}{$3.409$}};
					
					\node at (1,0.25) {\scalebox{0.8}{$0$}};
					\node at (4,0.4) {\scalebox{0.8}{$0.146$}};
					\node at (7,0.8) {\scalebox{0.8}{$0.589$}};
					\node at (10,0.4) {\scalebox{0.8}{$0.134$}};
					\node at (13,0.5) {\scalebox{0.8}{$0.252$}};
					
					\draw[ForestGreen,very thick] (1,0.897) -- (4,3.947) -- (7,4.007) -- (10,4.221) -- (13,3.409);
					\draw[dr,very thick] (1,0) -- (4,0.146) -- (7,0.589) -- (10,0.134) -- (13,0.252);
					
					\draw[black] (0,-1) -- (14,-1) -- (14,5) -- (0,5) -- (0,-1);
					
					\node at (1,-1.5) {Baboon};
					\node at (4,-1.5) {Malawi};
					\node at (7,-1.5) {Conference};
					\node at (10,-1.5) {Workplace};
					\node at (13,-1.5) {High school};
					
					\filldraw[ForestGreen] (1,-2.5) circle (3pt) {};
					\node at (5.8,-2.5) {proportion of barcodes for $H^{\inf}_1$ ending at infinity};
					\filldraw[dr] (1,-3) circle (3pt) {};
					\node at (6.67,-3) {proportion of additional barcodes for $\widehat{H}_1$ ending at infinity};
			\end{tikzpicture}}
		\end{center}
		\caption{Proportion of dimension 1 barcodes ending at infinity.}
	\end{figure}
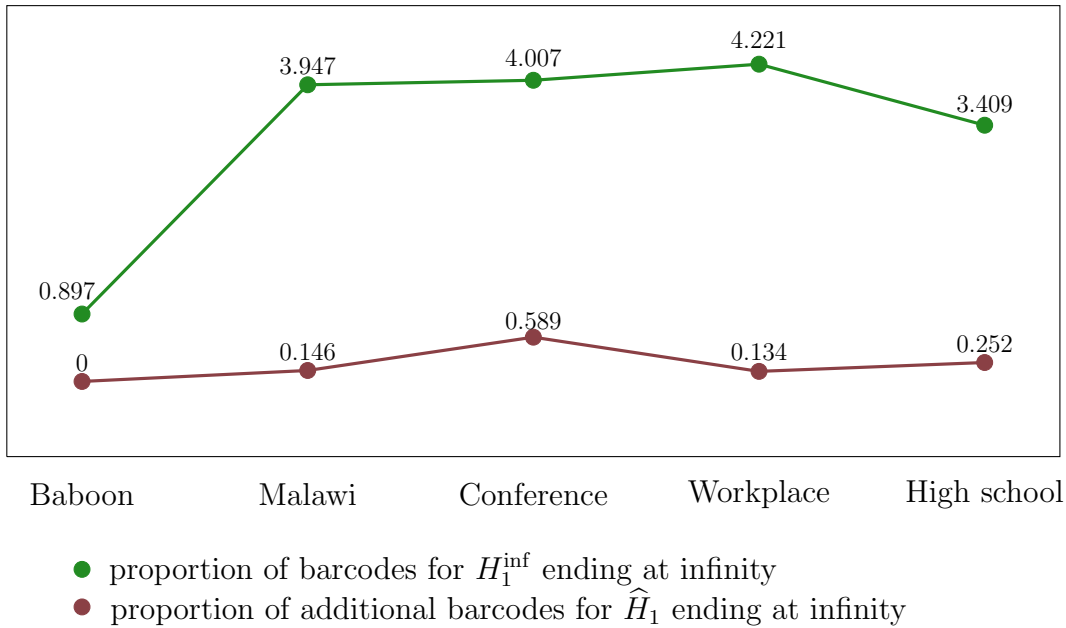
	
	For the data set of baboons, most persistent barcodes end before infinity. This result describes that baboons like to meet in groups with other individuals that they usually contact with. On the other hand, there is no $\hat{H}$ barcodes, showing that baboons never attend any group meetings if there is someone they do not and will not contact inside the group. All other data sets consist of human beings, and they all have relatively high proportion of persistent barcodes ending at infinity, meaning most pair-wise interactions do not end up with group discussion. The data set of conference holds the highest proportion of $\hat{H}$ barcodes ending at infinity, followed by the data set of high school. This phenomenon can be explained by that collaborations of more than 2 scientists or high school students are more common than those in rural Malawi and workplace. Moreover, such collaborations often do not require that individuals are familiar with each other in the sense of private connections.
	
	In the end of this section, we want to point out that there are many other potential statistics one can derive from persistent hypergraph model. For example, the lengths of dimension 1 barcodes reflect the differences in length between group meetings and pair-wise private meetings. The persistent betti numbers count the number of holes or anti-holes with each particular period length. Furthermore, higher dimensional barcodes reveal the contact information of group of more than 3 individuals. These statistic figures could potentially help understanding more deeply the behavior of face-to-face interactions under different environments.

	\section{Conclusion}
	
	The main contribution of this article is the introduction of persistent hypergraph model in data analysis. We introduce a new homology $\hat{H}$, and we extend classical algorithm of persistent homology to hypergraphs that computes both persistent $H^{\inf}$ and $\hat{H}$ barcodes simultaneously and efficiently. The persistent hypergraph model fits particularly well in the data sets of interaction networks. Both $H^{\inf}$ and $\hat{H}$ homology have real-world meanings, they represent connected components in the connection graph, cycles of pair-wise interactions and group meetings of at least 3 individuals. 
	
	The commonly used statistics features on face-to-face interaction analysis include connection graphs, mean values, standard deviations, distributions and their modifications. These features are essentially 0 or 1 dimensional in the topological viewpoint. Topological data analysis provides a very efficient way to reveal features in not only 0 and 1, but also higher dimensions. As demonstrated in the experiment, persistent hypergraph model reveals connected components, cycles, group meetings and other new and potentially useful information for understanding contact patterns.
	
	Furthermore, our algorithm for persistent hypergraph model gives both global features and local details, including information of the exact locations and the sizes of the cycles. Therefore, topological data analysis can be used as an important complementary to classical data analysis methods. Ideally, global and local results can be further combined with such as principal component analysis and machine learning to obtain a more complete and precise picture of the data sets.

\end{document}